\documentclass[12pt]{article}

\usepackage{amsmath, amsfonts, amssymb}
\usepackage{latexsym}

\newtheorem{thm}{Theorem}[section]
 \newtheorem{cor}[thm]{Corollary}

 \newtheorem{rem}[thm]{Remark}

\newcommand{\iy}{\infty}
\newcommand{\bC}{{\bf C}}
\newcommand{\bT}{{\bf T}}
\newcommand{\bZ}{{\bf Z}}
\newcommand{\bN}{{\bf N}}
\newcommand{\De}{\Delta}
\newcommand{\ka}{\kappa}
\newcommand{\ga}{\gamma}
\newcommand{\be}{\beta}
\newcommand{\ph}{\varphi}
\newcommand{\n}{\|}
\newcommand{\vsk}{\vspace{1mm}}

\newcommand{\ti}{\widetilde}

\begin{document}

\vspace*{1.5cm}
\begin{center}
{\Large\bf Szeg\"o via Jacobi}
\end{center}

\renewcommand{\thefootnote}{\fnsymbol{footnote}}

\vspace{2mm}
\begin{center}
{\bf A. B\"ottcher and H. Widom}\footnote[1]{The research of this author was supported by
National Science Foundation grant DMS-0243982.}
\end{center}

\vspace{1mm}
\begin{center}
{\em For Bernd Silbermann on His 65th Birthday}
\end{center}

\vspace{2mm}
\begin{quote}

\renewcommand{\baselinestretch}{1.0}
\footnotesize
{At present there exist numerous different approaches to results on
Toeplitz determinants of the type of Szeg\"o's strong limit theorem.
The intention of this paper is to show that Jacobi's theorem on the minors
of the inverse matrix remains one of the most comfortable tools for tackling
the matter. We repeat a known proof of the Borodin-Okounkov formula and
thus of the strong Szeg\"o limit theorem that is based on Jacobi's theorem.
We then use Jacobi's theorem to derive exact and asymptotic formulas for
Toeplitz determinants generated by functions with nonzero winding number.
This derivation is new and completely elementary.}
\end{quote}

\vspace{4mm}

\section{{\large Introduction}}

In \cite{CaPi3}, Carey and Pincus employ heavy machinery to establish a formula for Toeplitz
determinants generated by functions with nonvanishing winding number, and their paper begins
with the words ``Jacobi's theorem on the conjugate minors of the adjugate matrix formed from the
cofactors of the Toeplitz determinant has been the
main tool of previous attempts to generalize the classical strong Szeg\"o limit theorem.''
The purpose of the present paper is to demonstrate that Jacobi's theorem remains a perfect tool
for deriving Szeg\"o's theorem (which is no new message) and for generalizing the theorem
to the case of nonvanishing winding number (which seems to be not widely known).

\vsk
Let $\bT$ be the complex unit circle and let $f:\bT \to \bC\setminus\{0\}$ be a continuous function.
We denote $\ga \in \bZ$ the winding number of $f$ about the origin. So $f(t)=t^\ga a(t)$
($t \in \bT$) where $a$ has no zeros on $\bT$ and winding number zero. We define the
Fourier coefficients $f_k$ ($k \in \bZ$) of $f$ by
\[f_k=\frac{1}{2\pi}\int_0^{2\pi}f(e^{i\theta})e^{-ik\theta}d\theta\]
and consider the $n \times n$ Toeplitz matrices $T_n(f):=(f_{j-k})_{j,k=1}^n$
and their determinants $D_n(f):=\det T_n(f)$. We are interested in exact and asymptotic formulas
for $D_n(f)$.

\vsk
For the sake of definiteness, we assume that $a$ (equivalently, $f$) belongs to $C^\be$
with $\be > 1/2$, which means that $a$ has $[\be]$ continuous derivatives and that
the $[\be]$th derivative satisfies a H\"older condition with the exponent $\be-[\be]$.
To avoid well known subtleties, we suppose that $\be \notin \bN$.

\vsk
Under the above assumptions, $a$ has a logarithm $\log a$ in $C^\be$, and we denote the
Fourier coefficients of $\log a$ by $(\log a)_k$. We define $a_-$ and $a_+$ on $\bT$ by
\begin{equation}
a_-(t)=\exp \sum_{k=1}^\iy (\log a)_{-k}t^{-k}, \quad
a_+(t)=\exp \sum_{k=0}^\iy (\log a)_kt^k,\label{WHF}
\end{equation}
and we put $G(a):=\exp (\log a)_0$.
It is well known that $a_-^{\pm 1}$ and $a_+^{\pm 1}$ belong to $C^\beta$ together with $a$
(this results from the boundedness of the Cauchy singular integral operator on $C^\beta$
for $\be \notin \bN$). Clearly, $a=a_-a_+$.
This representation is called a Wiener-Hopf factorization of $a$.
The main actors in the following are the two functions
\begin{equation}
b=a_-a_+^{-1}, \quad c=a_+a_-^{-1}. \label{4a}
\end{equation}
Notice that $b \in C^\beta$, $c \in C^\beta$, and $bc=1$.

\vsk
For a continuous function $\ph$ on $\bT$, we define the infinite Toeplitz matrix $T(\ph)$ and the
infinite Hankel matrix $H(\ph)$ by $T(\ph):=(\ph_{j-k})_{j,k=1}^\iy$ and
$H(\ph):=(\ph_{j+k-1})_{j,k=1}^\iy$. These two matrices induce bounded linear operators on
$\ell^2({\bf N})$ whose (operator) norms satisfy $\n T(\ph) \n =\n \ph \n_\iy$ and
$\n H(\ph) \n \le \n \ph \n_\iy$, where $\n \cdot \n_\iy$ is the norm in $L^\iy(\bT)$.
We also define $\ti{\ph}$ by $\ti{\ph}\,(t)=\ph(1/t)$ for $t \in \bT$. Then
$H(\ti{\ph}\,)=(\ph_{-j-k+1})_{j,k=1}^\iy$. Let ${\cal P}_n$ be the linear space of
all trigonometric polynomials of degree at most $n$. If $\ph \in C^\be$, then there are
$p_n \in {\cal P}_n$ (the polynomials of best uniform approximation) such
that $\n\ph -p_n\n_\iy =O(n^{-\be})$. It follows that the $n$th singular number
$s_n$ of $H(\ph)$ satisfies
\[s_n\le \n H(\ph)-H(p_n)\n = \n H(\ph-p_n)\n \le \n \ph-p_n\n_\iy
=O(n^{-\be}),\]
which implies that $H(\ph)$ is a Hilbert-Schmidt operator for $\be > 1/2$. Consequently,
the product $H(b)H(\ti{c}\,)$ is a trace class operator and the determinant $\det (I-H(b)H(\ti{c}\,))$
is well-defined. We finally denote by $P_k$ and $Q_k$ the projections given by
\begin{eqnarray*}
& & P_k: (x_1, x_2, \ldots) \mapsto (x_1, \ldots, x_k, 0, 0, \ldots),\\
& & Q_k: (x_1, x_2, \ldots) \mapsto (0, \ldots, 0, x_{k+1}, x_{k+2}, \ldots).
\end{eqnarray*}

Here are the results we want to prove in this paper.

\begin{thm} \label{Th 1.1}
{\bf (Borodin-Okounkov formula)} The operator $I-H(b)H(\ti{c}\,)$
is invertible and
\[D_n(a)=G(a)^n\,\frac{\det (I-Q_nH(b)H(\ti{c}\,)Q_n)}{\det (I-H(b)H(\ti{c}\,))}\]
for all $n \ge 1$.
\end{thm}

\begin{thm} \label{Th 1.2}
{\bf (Szeg\"o's strong limit theorem)} We have
\[D_n(a)=G(a)^nE(a)(1+O(n^{1-2\be}))\]
where
\begin{eqnarray*}
E(a) & = & 1/\det(I-H(b)H(\ti{c}\,)) = \det T(a)T(a^{-1})\\
& = & \exp\sum_{k=1}^\iy k(\log a)_k(\log a)_{-k}
=  \exp\sum_{k=1}^\iy k(\log b)_k(\log c)_{-k}.
\end{eqnarray*}
\end{thm}

\begin{thm} \label{Th 1.3}
If $\ka >0$ and the matrix $T_{n+\ka}(a)$ is invertible, then the two operators
$I-H(b)H(\ti{c}\,)Q_{n+\ka}$ and $I-H(b)Q_nH(\ti{c}\,)Q_\ka$ are invertible and
\[D_n(t^{-\ka} a)=(-1)^{n\ka}D_{n+\ka}(a)F_{n,\ka}(a)\]
where
\begin{eqnarray}
F_{n,\ka}(a) & = & \det P_\ka T(t^{-n})\Big(I-H(b)H(\ti{c}\,)Q_{n+\ka}\Big)^{-1}T(b)P_\ka
\label{1.1}\\
& = & \det P_\ka \Big(I-H(b)Q_nH(\ti{c}\,)Q_{\ka}\Big)^{-1}T(t^{-n}b)P_\ka.
\label{1.2}
\end{eqnarray}
\end{thm}

\begin{thm} \label{Th 1.4}
{\bf (Fisher, Hartwig, Silbermann et al.)}
If $\ka >0$, then
\begin{equation}
F_{n,\ka}(a)= \det T_\ka(t^{-n}b)+O(n^{-3\be}) \label{1.3}
\end{equation}
and thus
\begin{equation}
D_n(t^{-\ka}a)=(-1)^{n\ka} G(a)^{n+\ka} E(a) \Big(\det T_\ka(t^{-n}b)+O(n^{-3\be})\Big)
\Big(1+O(n^{1-2\be})\Big). \label{1.4}
\end{equation}
\end{thm}

Proofs and comments on these theorems are in Sections 3, 4, 5, 6. In the following
Section~2 we recall Jacobi's theorem, and Section 7 contains additional material.

\section{{\large Jacobi's theorem}}

Let $C$ be an $m \times m$ matrix. For $i_1 < \ldots < i_s$ and $k_1 < \ldots < k_s$,
we denote by $C\left(\begin{array}{ccc}i_1 & \ldots & i_s\\k_1 & \ldots & k_s \end{array}\right)$
the determinant of the submatrix of $C$ that is formed by the intersection of the rows
$i_1, \ldots, i_s$ and the columns $k_1, \ldots, k_s$. We also define
the indices $j_1^\prime < \ldots <
j_{m-s}^\prime$ by $\{j_1^\prime, \ldots, j_{m-s}^\prime\}:=\{1, \ldots, m\} \setminus
\{j_1, \ldots, j_s\}$.

\begin{thm} \label{Th 2.1}
{\bf (Jacobi)} If $A$ is an invertible $m \times m$ matrix, then
\[A^{-1}\left(\begin{array}{ccc}i_1 & \ldots & i_s\\k_1 & \ldots & k_s \end{array}\right)
=(-1)^{\sum_{r=1}^s(i_r+k_r)}A\left(\begin{array}{ccc}k_1^\prime & \ldots & k_{m-s}^\prime\\
i_1^\prime & \ldots & i_{m-s}^\prime\end{array}\right) / \det A.\]
\end{thm}

A proof is in \cite{Gant}, for example. The following consequence of Theorem \ref{Th 2.1}
is from \cite{BoOk2}.

\begin{cor} \label{Cor 2.2}
If $K$ is a trace class operator and $I-K$ is invertible, then
\[\det P_n(I-K)^{-1}P_n = \frac{\det(I-Q_nKQ_n)}{\det(I-K)}\]
for all $n \ge 1$.
\end{cor}

{\em Proof.} If $m>n$ is sufficiently large, then $A:=I_{m \times m}-P_mKP_m$ is invertible together
with $I-K$. Theorem \ref{Th 2.1} with $\{i_1, \ldots, i_n\}=\{k_1, \ldots, k_n\} =\{1, \ldots, n\}$
applied to the $m \times m$ matrix $A$ yields
\[\det P_n(I_{m \times m}-P_mKP_m)^{-1}P_n
= \frac{\det(I_{(m-n)\times (m-n)}-Q_nP_mKP_mQ_n)}{\det(I_{m \times m}-P_mKP_m)},\]
which is equivalent to
\begin{equation}
\det P_n(I-P_mKP_m)^{-1}P_n
= \frac{\det(I-Q_nP_mKP_mQ_n)}{\det(I-P_mKP_m)}. \label{2.1}
\end{equation}
Since $P_mKP_m \to K$ in the trace norm as $ m \to \iy$ and the determinant is continuous on
identity minus trace class ideal, we may in (\ref{2.1}) pass to the limit $m \to \iy$ to get the desired formula.

\vsk
Here is another corollary of Jacobi's theorem. It was Fisher and Hartwig \cite{FiHa1}, \cite{FiHa2}
who were the first to write down this corollary and to recognize that it is the key to treating
the case of nonvanishing winding number.

\begin{cor} \label{Cor 2.3}
Let $\ka >0$ and suppose $T_{n+\ka}(a)$ is invertible. Then
\[D_n(t^{-\ka}a) =(-1)^{n\ka}D_{n+\ka}(a)\det (P_{n+\ka}-P_n)T_{n+\ka}^{-1}(a)P_\ka.\]
\end{cor}

{\em Proof.} This is immediate from Theorem \ref{Th 2.1} with $A=T_{n+\ka}(a)$
and
\begin{eqnarray*}
& & \det (P_{n+\ka}-P_n)T_{n+\ka}^{-1}(a)P_\ka = A^{-1}\left(\begin{array}{ccc}
n+1 & \ldots & n+\ka\\ 1 & \ldots & \ka\end{array}\right),\\
& & D_n(t^{-\ka}a)=A\left(\begin{array}{ccc}
\ka+1 & \ldots & \ka+n\\ 1 & \ldots & n\end{array}\right),\\
& & (-1)^{(n+\ka+1)\ka}=(-1)^{n\ka}.
\end{eqnarray*}

\section{{\large The Borodin-Okounkov formula}}

Theorem \ref{Th 1.1} was established by Borodin and Okounkov in \cite{BO}. Later it turned out
that (for positive functions $a$) it was already in Geronimo and Case's paper \cite{GeCa}.
The original proofs in \cite{BO}, \cite{GeCa} are quite complicated. Simpler proofs
were subsequently found in \cite{BaWi}, \cite{BoOk1}, \cite{BoOk2}. See also \cite{CaPi3}.
Here is the proof from \cite{BoOk2}, which is based on Jacobi's theorem.

\vsk
We apply Corollary \ref{Cor 2.2} to the trace class operator $K=H(b)H(\ti{c}\,)$.
The operator
\begin{eqnarray*}
I-K & = & T(bc)-H(b)H(\ti{c}\,) = T(b)T(c)\\
& = & T(a_-)T(a_+^{-1})T(a_-^{-1})T(a_+)
= T(a_-)T^{-1}(a_-a_+)T(a_+)
\end{eqnarray*}
has the inverse
\[(I-K)^{-1}=T(a_+^{-1}) T(a_-a_+)T(a_-^{-1})\]
and hence Corollary \ref{Cor 2.2} yields
\begin{equation}
\det P_n T(a_+^{-1}) T(a_-a_+)T(a_-^{-1})P_n
= \frac{\det (I-Q_nH(b)H(\ti{c}\,)Q_n)}{\det (I-H(b)H(\ti{c}\,))}. \label{3.1}
\end{equation}
Taking into account that
\[P_n T(a_+^{-1}) T(a_-a_+)T(a_-^{-1})P_n
= P_n T(a_+^{-1})P_n T(a_-a_+)P_nT(a_-^{-1})P_n\]
and that $P_n T(a_+^{-1})P_n$ and $P_n T(a_-^{-1})P_n$ are triangular with $1/G(a)$ and $1$,
respectively, on the main diagonal, we see that
the left-hand side of (\ref{3.1}) equals $D_n(a)/G(a)^n$.

\section{{\large The strong Szeg\"o limit theorem}}

We now prove Theorem \ref{Th 1.2}. As $H(b)H(\ti{c}\,)$ is in the trace class and
$Q_n=Q_n^*$ goes strongly to zero, we have $\det (I-Q_nH(b)H(\ti{c}\,)Q_n)
=1+o(1)$. Thus, Theorem \ref{Th 1.1} immediately gives
\[D_n(a)=G(a)^nE(a)(1+o(1))\quad \mbox{with} \quad
E(a)=1/\det(I-H(b)H(\ti{c}\,)). \]

To make the $o(1)$ precise, we proceed as in \cite{BoSi1}, \cite{BoSiBu}.
The $\ell$th singular number $s_\ell$ of $Q_nH(b)$ can be estimated by
\[s_\ell \le \n Q_nH(b)-Q_nH(p_{n+\ell})\n
=  \n Q_nH(b-p_{n+\ell})\n \le \n b-p_{n+\ell}\n_\iy\]
where $p_{n+\ell}$ is any polynomial in ${\cal P}_{n+\ell}$. There are such polynomials with
$\n b-p_{n+\ell}\n_\iy =O((n+\ell)^{-\be})$. This shows that the squared Hilbert-Schmidt norm
of $Q_nH(b)$ is
\[\sum_{\ell=0}^\iy s_\ell^2 =O\left(\sum_{\ell=0}^\iy (n+\ell)^{-2\be}\right)=O(n^{1-2\be}).\]
Thus, the Hilbert-Schmidt norm of $Q_nH(b)$ is $O(n^{1/2-\be})$. The same is true for the
operator $H(\ti{c}\,)Q_n$. Consequently, the trace norm of $Q_nH(b)H(\ti{c}\,)Q_n$
is $O(n^{1-2\be})$, which implies that
\[D_n(a)=G(a)^nE(a)(1+O(n^{1-2\be}))\quad \mbox{with} \quad
E(a)=1/\det(I-H(b)H(\ti{c}\,)). \]

We are left with the alternative expressions for $E(a)$. We start with
\begin{eqnarray}
& & 1/\det(I-H(b)H(\ti{c}\,))=1/\det T(b)T(c)=\det T^{-1}(c)T^{-1}(b)\nonumber \\
& & =\det T(a_+^{-1})T(a_-)T(a_+)T(a_-^{-1})
=\det T(a_-)T(a_+)T(a_-^{-1})T(a_+^{-1}). \label{4.1}
\end{eqnarray}
This equals
\[\det T(a_-a_+)T(a_-^{-1}a_+^{-1})=\det T(a)T(a^{-1}).\]
On the other hand, (\ref{4.1}) is
\begin{equation}
\det e^{T(\log a_-)} e^{T(\log a_+)} e^{-T(\log a_-)} e^{-T(\log a_+)} \label{4.2}
\end{equation}
and the Pincus-Helton-Howe formula \cite{HeHo}, \cite{Pi} (an easy proof of which was
recently found by Ehrhardt \cite{Ehr}) says that
\[\det e^A e^B e^{-A} e^{-B}=e^{{\rm tr}(AB-BA)}\]
whenever $A$ and $B$ are bounded and $AB-BA$ is in the trace class. Thus, (\ref{4.2})
becomes
\begin{eqnarray*}
& & \exp {\rm tr}\Big(T(\log a_-)T(\log a_+)-T(\log a_+)T(\log a_-)\Big)\\
& & = \exp {\rm tr} H(\log a_+)H((\log a_-)\,\ti{}\,)
=\exp \sum_{k=1}^\iy k(\log a_+)_k (\log a_-)_{-k}\\
& & = \exp  \sum_{k=1}^\iy k(\log a)_k (\log a)_{-k}
= \exp  \sum_{k=1}^\iy k(\log b)_k (\log c)_{-k}.
\end{eqnarray*}
Theorem \ref{Th 1.2} is completely proved.

\vsk
The treatment of the constant $E(a)$ given here is from \cite{WiII}.
For reviews of the gigantic development from Szeg\"o's original version of
his strong limit theorem \cite{Sz}
up to the present we refer to the books \cite{BoSiBu} and \cite{Simon}.

\section{{\large The exact formula for nonzero winding numbers}}

To prove Theorem \ref{Th 1.3} we use Corollary \ref{Cor 2.3} of Jacobi's theorem. Thus,
we must show that
\[\det (P_{n+\ka}-P_n)T_{n+\ka}^{-1}(a)P_\ka =: F_{n,\ka}(a)\]
is given by (\ref{1.1}) and (\ref{1.2}).

\vsk
We put $K:=H(b)H(\ti{c}\,)$, $m:=n+\ka$, $\De_n^\ka:=P_{n+\ka}-P_n$.
In \cite{BoOk1} it was shown
(in an elementary way) that the invertibility of $T_m(a)$ implies that $I-Q_mKQ_m$ is invertible and that
\[T_m^{-1}(a)=P_mT(a_+^{-1})\Big(I-T(c)Q_m(I-Q_mKQ_m)^{-1}Q_mT(b)\Big)T(a_-^{-1})P_m\]
(to get conformity with \cite{BoOk1} note that $P_mT(a_+^{-1})P_m= P_mT(a_+^{-1})$
and $P_mT(a_-^{-1})P_m=T(a_-^{-1})P_m$). We multiply this identity from the right by $P_\ka$
and from the left by $\De_n^\ka$. Since
\[T(a_-^{-1})P_mP_\ka=T(a_-^{-1})P_\ka =P_\ka T(a_-^{-1})P_\ka\]
and
\begin{eqnarray*}
& & \De_n^\ka P_mT(a_+^{-1})=\De_n^\ka T(a_+^{-1})= \De_n^\ka T(a_-^{-1})T(a_-)T(a_+^{-1})\\
& & = \De_n^\ka T(a_-^{-1})T(b)=\De_n^\ka T(a_-^{-1})\De_n^\ka T(b),
\end{eqnarray*}
we arrive at the formula
\begin{eqnarray*}
\det \De_n^\ka T_{m}^{-1}(a)P_\ka
&=&\det \De_n^\ka T(a_-^{-1})\De_n^\ka \cdot \det P_\ka T(a_-^{-1})P_\ka\\
& & \times \det \De_n^\ka T(b)\Big(I-T(c)Q_m(I-Q_mKQ_m)^{-1}Q_mT(b)\Big)P_\ka.
\end{eqnarray*}
As the matrix $T(a_-^{-1})$ is triangular with $1$ on the main diagonal,
we have
\begin{equation}
\det \De_n^\ka T(a_-^{-1})\De_n^\ka = \det P_\ka T(a_-^{-1})P_\ka =1 \label{5a}
\end{equation}
and are therefore left with the determinant of
\begin{equation}
\De_n^\ka T(b)\Big(I-T(c)Q_m(I-Q_mKQ_m)^{-1}Q_mT(b)\Big)P_\ka. \label{10}
\end{equation}
Taking into account that
$T(b)T(c)=T(bc)-H(b)H(\ti{c}\,)=I-K$
and $\De_n^\ka Q_m=0$, we obtain that (\ref{10}) equals
\begin{eqnarray}
& & \De_n^\ka T(b)P_\ka - \De_n^\ka(I-K)Q_m(I-Q_mKQ_m)^{-1}Q_mT(b)P_\ka\nonumber\\
& & = \De_n^\ka T(b) P_\ka + \De_n^\ka KQ_m(I-Q_mKQ_m)^{-1}Q_mT(b)P_\ka\nonumber\\
& & = \De_n^\ka \Big(I+KQ_m(I-Q_mKQ_m)^{-1}Q_m\Big)T(b)P_\ka. \label{11}
\end{eqnarray}
Since
\[(I-Q_mKQ_m)^{-1}Q_m (I-Q_mKQ_m)
=(I-Q_mKQ_m)^{-1}(I-Q_mKQ_m)Q_m=Q_m,\]
we get
$(I-Q_mKQ_m)^{-1}Q_m=Q_m(I-Q_mKQ_m)^{-1}$
and thus
\begin{equation}
I+KQ_m(I-Q_mKQ_m)^{-1}Q_m = I+KQ_m(I-Q_mKQ_m)^{-1}. \label{12a}
\end{equation}
We have $I-KQ_m=(I-Q_mKQ_m)(I-P_mKQ_m)$ and the operators $I-Q_mKQ_m$ and $I-P_mKQ_m$
are invertible; note that $(I-P_mKQ_m)^{-1}=I+P_mKQ_m$. Consequently, $I-KQ_m$ is also
invertible. It follows that (\ref{12a}) is
\begin{eqnarray}
& & (I-Q_mKQ_m)(I-Q_mKQ_m)^{-1}+KQ_m(I-Q_mKQ_m)^{-1}\nonumber\\
& & = (I-Q_mKQ_m+KQ_m)(I-Q_mKQ_m)^{-1} \nonumber\\
& & = (I+P_mKQ_m)(I-Q_mKQ_m)^{-1}\nonumber\\
& & = (I-P_mKQ_m)^{-1} (I-Q_mKQ_m)^{-1}= (I-KQ_m)^{-1}. \label{12}
\end{eqnarray}
In summary, (\ref{11}) is $\De_n^\ka(I-KQ_m)^{-1}T(b)P_\ka$ and we have proved the theorem
with
\[F(n,\ka)=\det\De_n^\ka (I-KQ_{n+\ka})^{-1}T(b)P_\ka.\]
The operator $T(t^{-k})$ sends $(x_1, x_2, \ldots)$ to $(x_{k+1}, x_{k+2}, \ldots)$.
It follows that $\De_n^\ka =P_\ka T(t^{-n})$, which yields (\ref{1.1}).

\vsk
The matrix $I-KQ_{n+\ka}$ is of the form
\[\left(\begin{array}{cc}I_{(n+\ka)\times (n+\ka)} & *\\0 & B\end{array}\right),\]
and since $I-KQ_{n+\ka}$ is invertible, the matrix $B$ must also be invertible.
The matrix
\begin{equation}
M:=T(t^{-n})(I-KQ_{n+\ka})T(t^n) \label{13}
\end{equation}
results from $I-KQ_{n+\ka}$ by deleting the first $n$ rows and first $n$ columns.
Consequently, $M$ has the form
\[\left(\begin{array}{cc}I_{\ka\times \ka} & *\\0 & B\end{array}\right),\]
and the invertibility of $B$ implies that $M$ is invertible.
Since $T(t^k)T(t^{-k})=Q_k$ and hence
\[ MT(t^{-n})=T(t^{-n})(I-KQ_{n+\ka})Q_n
= T(t^{-n})(I-KQ_{n+\ka}),\]
we get
$T(t^{-n})(I-KQ_{n+\ka})^{-1}=M^{-1}T(t^{-n})$.
Inserting this in (\ref{1.1}) we arrive at the formula
\[F(n,\ka)=\det P_\ka M^{-1}T(t^{-n})T(b)P_\ka =\det P_\ka M^{-1}T(t^{-n}b)P_\ka.\]
Finally, the identity $T(t^{-k})H(\varphi)=H(\varphi)T(t^k)$
shows that
\begin{eqnarray*}
M & = & I-T(t^{-n})H(b)H(\ti{c}\,)T(t^{n+\ka})T(t^{-n-\ka})T(t^n)\\
& = & I-T(t^{-n})H(b)H(\ti{c}\,)T(t^{n})T(t^{\ka})T(t^{-\ka})\\
& = & I-H(b)T(t^n)T(t^{-n})H(\ti{c}\,)T(t^{\ka})T(t^{-\ka})\\
& = & I-H(b)Q_nH(\ti{c}\,)Q_\ka,
\end{eqnarray*}
which gives (\ref{1.2}) and completes the proof of Theorem \ref{Th 1.3}.

\vsk
A result like Theorem \ref{Th 1.3} appeared probably first in \cite{OT48}.
Given a set $E \subset \bZ$, we denote by $P_E$ the projection on $L2(\bT)$
defined by
\[P_E: \sum_{k \in \bZ} x_k t^k \mapsto \sum_{k \in E} x_k t^k,\]
and for a function $\ph$ on $\bT$, we denote the operator of multiplication
by $\ph$ on $L2(\bT)$ also by $\ph$. Let $U$ and $V$ be the operators on $L2(\bT)$
given by
\[U:=P_{\{1,2,\ldots\}}t^{-n-\ka+1}b, \quad V:=P_{\{-1,-2,\ldots\}}t^{n+\ka-1}c.\]
One can show that $I-VU$ is invertible. Put $Y=(Y_{ij})_{i,j=0}^{\ka-1}$ with
\[Y_{ij}=P_{\{-i\}} t^{-n-\ka+1}b(I-VU)^{-1}P_{\{j\}}.\]
Lemma 3.2 of \cite{OT48} says that
\begin{equation}
\det \De_n^\ka T_{n+\ka}^{-1}(a)P_\ka =(-1)^{\ka}\det Y, \label{Harold}
\end{equation}
which together with Corollary \ref{Cor 2.3} yields
\[D_n(t^{-\ka}a)=(-1)^{n\ka}D_{n+\ka}(a)(-1)^{\ka} \det Y.\]
Clearly, this highly resembles Theorem \ref{Th 1.3}. In Remark \ref{Rem 6.1}
we will show that the right-hand side of (\ref{Harold}) indeed coincides with
(\ref{1.2}).

\vsk
Carey and Pincus \cite{CaPi3} state that
\begin{equation}
D_n(t^{-\ka}a)=(-1)^{n\ka}G(a)^{n+\ka}E(a)\ti{F}_{n,\ka}(a)(1+O(n^{1-2\be})) \label{W2}
\end{equation}
with
\begin{equation}
\ti{F}_{n,\ka}(a)=\det P_\ka\Big(I-H(b)Q_{n-\ka}H(\ti{c}\,)\Big)^{-1}T(t^{-n}b)P_\ka. \label{W3}
\end{equation}
The proof of (\ref{W2}), (\ref{W3}) given in \cite{CaPi3} is complicated and based on the methods
developed in these authors' work \cite{CaPi1}, \cite{CaPi2}, \cite{CaPi3}.
We will return to
(\ref{W2}), (\ref{W3}) in Remark \ref{Rem 6.6}.

\section{{\large The asymptotic formula for nonzero winding numbers}}

Theorem \ref{Th 1.4} is an easy consequence of Theorem \ref{Th 1.3}. Since
$\n H(b)Q_nH(\ti{c}\,)Q_\ka \n \to 0$ as $n \to \iy$, we obtain that
\begin{equation}
F_{n,\ka}(a)=\det\Big[P_\ka T(t^{-n}b)P_\ka
+\sum_{k=1}^\iy P_\ka \Big(H(b)Q_nH(\ti{c}\,)Q_\ka\Big)^kT(t^{-n}b)P_\ka\Big]
\label{W4}
\end{equation}
for all sufficiently large $n$. We know that there are polynomials $p_n$ and $q_n$ in
${\cal P}_{n-\ka-1}$ such that $\n b-p_n\n_\iy =O(n^{-\be})$ and
$\n c-q_n\n_\iy =O(n^{-\be})$. It follows that
\begin{eqnarray*}
& & \n H(b)Q_n\n =\n H(b-p_n)Q_n\n \le \n b-p_n\n_\iy =O(n^{-\be}),\\
& & \n Q_nH(\ti{c}\,)\n =\n Q_n H(\ti{c}-q_n)\n \le \n c-q_n\n_\iy =O(n^{-\be}),\\
& & \n T(t^{-n}b)P_\ka\n =\n T(t^{-n}(b-p_n))P_\ka\n = \n b-p_n\n_\iy =O(n^{-\be}).
\end{eqnarray*}
Since $P_\ka$ is a trace class operator, the sum in (\ref{W4}) is $O(n^{-3\be})$
in the trace norm, which implies the claim of Theorem \ref{Th 1.4}.

\vsk
We remark that a result close to Theorem \ref{Th 1.4} was already established by Fisher and Hartwig
\cite{FiHa1}, \cite{FiHa2} using different methods. Theorem \ref{Th 1.4} as it is stated,
a formula similar to (\ref{W4}), and the estimates via $\n b-p_n\n_\iy$ and
$\n c-q_n\n_\iy$ used above are due to Silbermann and one of the authors \cite{BoSi1}.

\section{{\large Remarks}}

Here are a few additional issues.

\begin{rem} \label{Rem 6.2}
{\rm The proof of Theorem \ref{Th 1.3} given in Section 5
was done under the minimal assumption that $T_{n+\ka}(a)$
be invertible.
The proof can be simplified if one is satisfied by the formula for
sufficiently large $n$ only. Indeed, the operator $K=H(b)H(\ti{c}\,)$ is compact and hence
$\| KQ_m\| \to 0$ as $m \to \iy$. It follows that $\| KQ_m \| <1$ whenever $m=n+\ka$ is large enough,
and for these $m$ we can replace all between (\ref{11}) and (\ref{12}) by the simple series argument
\begin{eqnarray*}
& & I+KQ_m(I-Q_mKQ_m)^{-1}Q_m\\
& & = I +KQ_m + KQ_m KQ_m +KQ_m KQ_m KQ_m + \ldots \\
& & =(I-KQ_m)^{-1}.
\end{eqnarray*}
Moreover, if $\| KQ_m \| <1$ then the invertibility of the operator
(\ref{13}) is obvious and we can omit the piece of the proof dedicated to the invertibility of
(\ref{13}).}
\end{rem}

\begin{rem} \label{Rem 6.1}
{\rm We prove that the right-hand side of (\ref{Harold}) is the same as (\ref{1.2}).
We identify $L2(\bT)$ with $\ell^2(\bZ)$ in the natural fashion and think of operators
on $L2(\bT)$ as acting by infinite matrices on $\ell^2(\bZ)$. Let $m:=n+\ka$.
The matrices of $U$ and $V$
have the entries
\[U_{ij}=\left\{\begin{array}{lll}b_{i-j+m-1} & \mbox{if} & i>0,\\
0 & \mbox{if} & i \le 0,\end{array}\right.\quad
V_{ij}=\left\{\begin{array}{lll}c_{i-j-m+1} & \mbox{if} & i<0,\\
0 & \mbox{if} & i \ge 0,\end{array}\right.\]
and the matrix of the multiplication operator  $B:=t^{-m+1}b$ has
$i,j$ entry $b_{i-j+m-1}$. The $i,j$ entry of the product $VU$
equals
\[\sum_{k>0} \ti{c}_{-i+k+m-1}b_{k-j+m-1}\]
for $i<0$ and is $0$ for $i \ge 0$. If we set
\[H_{ij}=(VU)_{ij} \quad (j<0), \quad L_{ij}=(VU)_{ij} \quad (j \ge 0),\]
with both equal to $0$ when $i \ge 0$, then the operator $I-VU$ has the
matrix representation
\begin{equation}
\left(\begin{array}{cc} I-H & -L\\ 0 & I \end{array}\right) \label{W1}
\end{equation}
corresponding to the decomposition $\ell^2(\bZ)=\ell^2(\bZ_-) \oplus \ell^2(\bZ_+)$
with $\bZ_-=\{-1,-2, \ldots\}$ and $\bZ_+=\{0,1,2, \ldots\}$. The inverse of (\ref{W1})
is
\[\left(\begin{array}{cc} (I-H)^{-1} & (I-H)^{-1}L\\ 0 & I \end{array}\right).\]
Hence the $i,j$ entry of $Y$ is
\[Y_{ij} =(BP_-(I-H)^{-1}L)_{-i,j}+b_{-i-j+m-1},\]
where $P_\pm:=P_{\bZ_\pm}$. Here $i$ and $j$ run from $0$ to $\ka-1$. Now replace $i$ by
$\ka-i-1$. The new index also runs from $0$ to $\ka-1$. Thus,
\[Y_{\ka-i-1,j} =(BP_-(I-H)^{-1}L)_{-\ka+i+1,j}+b_{-\ka+i-j+m}.\]
Let $J$ be given on $\ell^2(\bZ)$ by $(Jx)_k=x_{k-1}$. Then $J2=I$ and $P_-J=JP_+$.
Consequently,
\begin{eqnarray*}
Y_{\ka-i-1,j} & = & (BP_-J(I-JHJ)^{-1}JL)_{-\ka+i+1,j}+b_{-\ka+i-j+m}\\
& = & (BJP_+(I-JHJ)^{-1}JL)_{-\ka+i+1,j}+b_{-\ka+i-j+m}.
\end{eqnarray*}
The matrix $(b_{-\ka+i-j+m})$ at the end delivers $T(t^{\ka-m}b)$. Next, $BJ$
has $-\ka+i+1,j$ entry $b_{i+j+m-\ka+1}$, which is the $i,j$ entry of
$H(t^{\ka-m}b)$. The $i,j$ entry of $JHJ$ is
\[\sum_{k \ge 0}\ti{c}_{i+k+m+1}b_{k+j+m+1}
=\sum_{k \ge \ka} \ti{c}_{i+k+m-\ka+1}b_{k+j+m-\ka+1}\]
and so the operator itself is $H(t^{\ka-m}\ti{c}\,)Q_\ka H(t^{\ka-m}b)$. Finally,
the $i,j$ entry of $JL$ is equal to
\[\sum_{k \ge 0} \ti{c}_{i+k+m+1}b_{k-j+m}=\sum_{k \ge \ka}
\ti{c}_{i+k-\ka+m+1}b_{k-\ka-j+m}\]
whence $JL=H(t^{\ka-m}\ti{c}\,)Q_\ka T(t^{\ka-m}b)$. Let $D:=H(t^{\ka-m}b)$ and
$C:=H(t^{\ka-m}\ti{c}\,)$. We have shown that $Y_{\ka-i-1,j}$ is the $i,j$ entry of
\begin{eqnarray*}
& & \Big(D(I-CQ_\ka D)^{-1}CQ_\ka+I\Big)T(t^{\ka-m}b)\\
& & = (I-DCQ_\ka)^{-1}T(t^{\ka-m}b)\\
& & =\Big(I-H(t^{\ka-m}b)H(t^{\ka-m}\ti{c}\,)Q_\ka\Big)^{-1}T(t^{\ka-m}b)\\
& & =\Big(I-H(b)Q_{m-\ka}H(\ti{c}\,)Q_\ka\Big)^{-1}T(t^{\ka-m}b)\\
& & =\Big(I-H(b)Q_{n}H(\ti{c}\,)Q_\ka\Big)^{-1}T(t^{-n}b).
\end{eqnarray*}
It follows that $(-1)^{\ka}\det(Y_{ij})=\det(Y_{\ka-i-1,j})$ equals
\[\det P_\ka \Big(I-H(b)Q_nH(\ti{c}\,)Q_\ka\Big)^{-1}T(t^{-n}b)P_\ka,\]
as desired.}
\end{rem}

\begin{rem} \label{Rem 6.6}
{\rm We show that (\ref{W2}), (\ref{W3}) are consistent with Theorem \ref{Th 1.3}.
Let first
\[M_n=I-H(b)Q_nH(\ti{c}\,)Q_\ka, \quad R_n=I-H(b)Q_nH(\ti{c}\,), \quad T_n=T(t^{-n}b).\]
We have $M_n=R_n+H(b)Q_nH(\ti{c}\,)P_\ka=:R_n+Z_n$. Using best approximation of $b$ and $c$
as above, we get $\n Z_n \n =O(n^{-2\be})$, and it is clear that $\n R_n^{-1}\n =O(1)$.
The identity $M_n^{-1}=(I+R_n^{-1}Z_n)^{-1}R_n^{-1}$ implies that
\[P_\ka M_n^{-1}T_nP_\ka=P_\ka(I+R_n^{-1}Z_n)^{-1}P_\ka R_n^{-1}T_nP_\ka
+P_\ka(I+R_n^{-1}Z_n)^{-1}Q_\ka R_n^{-1}T_nP_\ka,\]
and the second term on the right is zero because $P_\ka(I+R_n^{-1}Z_n)^{-1}$
has $P_\ka$ at the end.
It follows that
\[\det P_\ka M_n^{-1}T_nP_\ka = \det P_\ka R_n^{-1}T_nP_\ka \,(1+O(n^{-2\be})),\]
or equivalently,
\begin{equation}
F_{n,\ka}(a)=\det P_\ka \Big(I-H(b)Q_nH(\ti{c}\,)\Big)^{-1}T(t^{-n}b)P_\ka \,(1+O(n^{-2\be})).
\label{Con1}
\end{equation}
To change the $Q_n$ to $Q_{n-\ka}$, let
\[S_n=I-H(b)Q_{n-\ka}H(\ti{c}\,), \quad X_n=H(b)(Q_{n-\ka}-Q_n)H(\ti{c}\,).\]
Then $R_n=S_n+X_n$, $\n X_n\n =O(n^{-2\be})$, $\n S_n^{-1}\n =O(1)$, and
\[P_\ka R_n^{-1}T_nP_\ka=P_\ka(I+S_n^{-1}X_n)^{-1}P_\ka S_n^{-1}T_nP_\ka
+P_\ka(I+S_n^{-1}X_n)^{-1}Q_\ka S_n^{-1}T_nP_\ka.\]
This time the second term on the right does not disappear and hence all we can say
is that
\[\det P_\ka R_n^{-1}T_nP_\ka = \det P_\ka S_n^{-1}T_nP_\ka + O(n^{-2\be}).\]
Combining this and (\ref{Con1}) we arrive at the formula
$F_{n,\ka}(a)=\ti{F}_{n,\ka}(a)+O(n^{-2\be})$
and thus at
\[D_n(t^{-\ka} a)=(-1)^{n\ka}G(a)^{n+\ka}E(a)(\ti{F}_{n,\ka}(a)+O(n^{-2\be}))(1+O(n^{1-2\be})),\]
which is not yet (\ref{W2}) but reveals that (\ref{W2}) is consistent with Theorem \ref{Th 1.3}.
We emphasize that (\ref{W2}) is an asymptotic result while Theorem \ref{Th 1.3} provides us
with an exact formula. Moreover, $F_{n,\ka}(a)$ is a little better than $\ti{F}_{n,\ka}(a)$
since $Q_n$ and $Q_\ka$ are ``smaller'' than $Q_{n-\ka}$ and $I$.}
\end{rem}

\begin{rem} \label{Rem 6.3}
{\rm
We worked with the Wiener-Hopf factorization $a=a_-a_+$ specified by (\ref{WHF}).
One can do everything
if one starts with an arbitrary Wiener-Hopf
factorization $a=a_-a_+$. The different factorizations are all of the form
$a=(\mu^{-1}a_-)(\mu a_+)$ where $\mu$ is a nonzero complex number. The functions
$b$ and $c$ are then defined by
\[b=(\mu^{-1}a_-)(\mu a_+)^{-1}=\mu^{-2}a_-a_+^{-1},
\quad c=(\mu^{-1}a_-)^{-1}(\mu a_+)=\mu2 a_-^{-1}a_+.\]
Theorems \ref{Th 1.1} and \ref{Th 1.2} are invariant under this change. The only difference in
Theorems \ref{Th 1.3} and \ref{Th 1.4} is that if we replace $a_-$ by $\mu^{-1}a_-$ in
(\ref{5a}), then the determinants are $\mu^\ka$ and their product becomes $\mu^{2\ka}$.
Since $G(a)=G(a_-a_+)$ and $G(c)=\mu^2G(a_-a_+)$, we obtain that $\mu2=G(a)^{-\ka}G(c)^\ka$
and hence
\[\det \De_n^\ka T_{n-\ka}^{-1}(a)P_\ka =G(a)^{-\ka}G(c)^\ka F_{n,\ka}(a).\]
The invariant versions of Theorems \ref{Th 1.3} and \ref{Th 1.4} are
\begin{eqnarray*}
D_n(t^{-\ka}a)
& = & (-1)^{n\ka}D_{n+\ka}(a)G(a)^{-\ka} G(c)^{\ka}F_{n,\ka}(a)\\
& = & (-1)^{n\ka}G(a)^n E(a) G(c)^{\ka}\Big(\det T_\ka(t^{-n}b)+O(n^{-3\be})\Big)
(1+O(n^{1-2\be})),
\end{eqnarray*}
where $F_{n,k}(a)$ is given by (\ref{1.1}) and (\ref{1.2}).}
\end{rem}

\begin{rem} \label{Rem 6.4}
{\rm
Theorems \ref{Th 1.1} to \ref{Th 1.4} can be extended to block
Toeplitz operators generated by $\bC^{N \times N}$-valued $C^\beta$-functions.
In that case one has to start with two Wiener-Hopf factorizations $a=u_-u_+=v_+v_-$ and to put
$b=v_-u_+^{-1}$, $c=u_-^{-1}v_+$. Theorem \ref{Th 1.1} and its proof remain in force literally.
Theorem \ref{Th 1.2} and its proof yield the operator determinants for $E(a)$ but
not the expressions in terms of the Fourier coefficients of $\log a$, $\log b$, $\log c$.
In Theorems \ref{Th 1.3} and \ref{Th 1.4} one has to require that all partial indices
be equal to one another. The result reads
\begin{eqnarray*}
& & D_n\left[\left(\begin{array}{ccc}
t^{-\ka} & &\\& \ddots & \\ & & t^{-\ka} \end{array}\right)a\right]\\
& & = (-1)^{n\ka N} D_{n+\ka}(a)G(a)^{-\ka}G(c)^\ka
\det P_\ka \Big(I-H(b)Q_nH(\ti{c}\,)Q_\ka\Big)^{-1}T(t^{-n}b)P_\ka\\
& & = (-1)^{n\ka N} G(a)^{n}E(a)G(c)^\ka
\Big(\det T_\ka(t^{-n}b)+O(n^{-3\be})\Big) (1+O(n^{1-2\be})).
\end{eqnarray*}
For details see \cite{BoOk1} and \cite{BoSiBu}.}
\end{rem}

\begin{rem} \label{Rem 6.5}
{\rm The case of positive winding numbers can be reduced to negative winding numbers
by passage to transposed matrices because
$D_n(t^\ka a)=D_n(t^{-\ka}\ti{a}\,)$. Let $a=a_-a_+$ be any Wiener-Hopf factorization.
We denote the functions associated with $\ti{a}$ through (\ref{4a}) by
$b_*$ and $c_*$:
\[b_*=\ti{a_+}\,\ti{a_-}^{-1}=\ti{c}\,, \quad c_*=\ti{a_-}\,\ti{a_+}^{-1}=\ti{b}.\]
{From} Remark \ref{Rem 6.3} we infer that
\[D_n(t^\ka a)=(-1)^{n\ka} D_{n+\ka}(a)G(a)^{-\ka}G(c_*)^\ka F_{n,\ka}(\ti{a})\]
with $G(c_*)=G(\ti{b})=G(b)$ and
\begin{eqnarray*}
F_{n,\ka}(\ti{a}) & = & \det P_\ka \big(I-H(b_*)Q_nH(\ti{c_*})Q_\ka\Big)^{-1}T(t^{-n}b_*)P_\ka\\
& = & \det P_\ka \Big(I-H(\ti{c}\,)Q_nH(b)Q_\ka\Big)^{-1}T(t^{-n}\ti{c}\,)P_\ka\\
& = & \det P_\ka T(t^nc)\Big(I-Q_\ka H(b)Q_nH(\ti{c}\,)\Big)^{-1}P_\ka.
\end{eqnarray*}
}
\end{rem}

\vspace{4mm}
\noindent
\begin{minipage}[t]{7cm}
A. B\"ottcher\\
Fakult\" at f\" ur Mathematik \\
Technische Universit\"at Chemnitz \\09107 Chemnitz\\ Germany \\[0.5ex]
aboettch@mathematik.tu-chemnitz.de
\end{minipage}
\hspace{2cm}\begin{minipage}[t]{6cm}
H. Widom\\
Department of Mathematics\\
University of California\\
Santa Cruz, CA 95064\\
USA\\[0.5ex]
widom@math.ucsc.edu
\end{minipage}

\vspace{8mm}
\noindent
MSC 2000: 47B35

\end{document}